\documentclass[10pt]{amsart}
\usepackage{amsmath,amscd,amssymb,epsfig}

\newcommand{\la}{{\lambda}}
\newcommand{\hf}{{\hat{f}}}
\newcommand{\hg}{{\hat{g}}}

\newtheorem{theorem}{Theorem}[section]
\newtheorem{corollary}[theorem]{Corollary}
\newtheorem{proposition}[theorem]{Proposition}
\newtheorem{conjecture}[theorem]{Conjecture}
\newtheorem{lemma}[theorem]{Lemma}

\newtheorem{example}[theorem]{Example}

\newtheorem{question}[theorem]{Question}

\newtheorem{observation}[theorem]{Observation}

\numberwithin{figure}{section}

\def\mystrut{\hbox{\vrule height8.0pt depth2pt width0pt}}
\def\mybox{\hbox to 10.0pt}
\def\norulefill{\leaders\hrule height0pt\hfill}
\def\nr#1{\multispan{#1}\norulefill}
\def\hr#1{\multispan{#1}\hrulefill}

\begin{document}
\title[SYT of Truncated Shapes]%
{Enumeration of standard Young tableaux\\ of certain truncated
shapes \footnote{\, \\Dedicated to Doron Zeilberger on the
occasion of his ${\rm 60^{th}}$ birthday. {\em Mazal Tov!}}}

\author{Ron M.\ Adin}
\address{Department of Mathematics\\
Bar-Ilan University\\
Ramat-Gan 52900\\
Israel} \email{radin@math.biu.ac.il}

\author{Ronald C.\ King}
\address{School of Mathematics\\
University of Southampton\\
Southampton SO17 1BJ\\
England} \email{R.C.King@soton.ac.uk}

\author{Yuval Roichman}
\address{Department of Mathematics\\
Bar-Ilan University\\
Ramat-Gan 52900\\
Israel} \email{yuvalr@math.biu.ac.il}

\thanks{Two of the authors, RMA and YR, were partially supported by 
internal research grants from Bar-Ilan University.}

\date{Submitted: December 16, 2010;\ Revised: August 16, 2011.}

\begin{abstract}
Unexpected product formulas for the number of standard Young
tableaux of certain truncated shapes are found and proved. 
These include
shifted staircase shapes minus a square in the NE corner,
rectangular shapes minus a square in the NE corner, and some
variations.
\end{abstract}

\maketitle

\section{Introduction}

A truncated shape is obtained from 
a Ferrers diagram
(in the English notation, where parts decrease from top to bottom)
by deleting cells from the NE corner.
Interest in the enumeration of standard Young
tableaux of truncated shapes is enhanced by a recent
result~\cite[Prop. 9.7]{AR_tft2}: the number of geodesics between
distinguished pairs of antipodes in the flip graph of
triangle-free triangulations is equal to twice the number of Young
tableaux of a truncated shifted staircase shape. Motivated by this
result, 
extensive computations were carried out for the 
number of standard Young tableaux of these and other truncated
shapes. It was found that, in certain distinguished cases, all
prime factors of these numbers are relatively small, hinting at the existence
of product formulas. In this paper, 
such product formulas  are proved for rectangular and shifted staircase shapes
truncated by a square, or nearly a square.

\smallskip

A different method to derive product formulas, for other families
of truncated shapes, has independently been developed by G.\ Panova~\cite{Panova}.

\smallskip

The rest of the paper is organized as follows. Basic concepts are
described in Section~\ref{sec:basic}. The general idea of pivoting
is presented in Section~\ref{sec:idea}.
Section~\ref{sec:staircase} contains detailed proofs for truncated
shifted staircase shapes, with Theorem~\ref{main1111} as the main
result; while Section~\ref{sec:rectangular} contains an analogous
development for truncated rectangular shapes, with
Theorem~\ref{main2222} as the main result. Section~\ref{sec:final}
contains final remarks and open problems.

\section{Preliminaries and Basic Concepts}\label{sec:basic}

A {\em partition} $\lambda$ of a positive integer $N$ is a sequence of
non-negative integers $(\lambda_1,\lambda_2,\ldots, \lambda_m)$
such that 
$\lambda_1\ge \lambda_2\ge \ldots\ge \lambda_m\ge 0$
and the total size 
$|\la|:=\lambda_1+\lambda_2+\ldots+\lambda_m$ is $N$.
The {\em Ferrers diagram} $[\lambda]$ of shape $\lambda$ is a left-justified array 
of $N$ cells, with row $i$ (from top to bottom) containing $\lambda_i$ cells. 
A {\em standard Young tableau} (SYT) $T$ of shape $\lambda$ is a labeling by
$\{1,2,\ldots,N\}$ of the cells in the diagram $[\lambda]$ such that
every row is increasing from left to right, and every
column is increasing from top to bottom. 
The number of SYT of shape $\lambda$ is denoted by $f^\lambda$.

\begin{proposition}\label{t.num_regular_SYT}{\rm\ (The Frobenius-Young Formula)}~\cite{F, Y}
The number of SYT of shape $\la=(\lambda_1,\lambda_2,\ldots,\lambda_m)$ 
with $\lambda_1\geq \lambda_2\geq \ldots\geq \lambda_m\geq 0$ is
$$
f^{\la} = \frac{(|\la|)!}{\prod_{i=1}^{m} (\la_i+m-i)!} \cdot
\prod_{1\le i<j \le m} (\la_i - \la_j - i + j).
$$
\end{proposition}

Note that adding trailing zeros to $\la$ (with $m$ appropriately
increased) does not affect the right-hand side of the above formula.

\bigskip

A partition $(\lambda_1,\lambda_2,\ldots, \lambda_m)$ of $N$ is
{\em strict} if $\lambda _1>\lambda_2>\ldots>\lambda_m >0$.
The corresponding diagram of {\em shifted shape} $\lambda$
is the array of $N$ cells with row $i$ containing $\lambda_i$ cells
and indented $i-1$ places. 
A standard Young tableau (SYT) $T$ of shifted shape $\lambda$ is 
a labeling by $\{1,2,\ldots,N\}$ of the cells in the
diagram $[\lambda]$ such that each row and column is increasing. 
The number of SYT of shifted shape $\lambda$ is denoted by $g^\lambda$.

\begin{proposition}\label{t.num_shifted_SYT}{\rm\ (Schur's
Formula)}\cite{Schur}\cite[p. 267 (2)]{Md}
The number of SYT of shifted shape $\la=(\lambda_1,\lambda_2,\ldots,\lambda_m)$ with
$\lambda_1> \lambda_2> \ldots>\lambda_m>0$ is
$$
g^{\la} = \frac{(|\la|)!}{\prod\limits_{i=1}^m \la_i!} \cdot
\prod_{1\le i<j\le m} \frac{\la_i - \la_j}{\la_i + \la_j}.
$$
\end{proposition}



\medskip

It is well known that both $f^\la$ and $g^\la$ also have hook
length formulas, but the equivalent formulas above will be more
convenient for our calculations.

\medskip

We shall frequently use the following two basic operations on partitions.
The {\em union} $\la \cup \mu$ of two partitions $\la$ and $\mu$ is 
simply their multiset union. 
We shall usually assume that 
each part of $\la$ is greater or equal than each part of $\mu$.
The {\em sum} of two partitions $\la=(\la_1,\ldots,\la_m)$ and $\mu=(\mu_1,\ldots,\mu_m)$
(with trailing zeros added in order to get the same number of parts) is
$$
\la + \mu := (\la_1+\mu_1, \ldots, \la_m+\mu_m).
$$

\medskip

For any nonnegative integer $m$, let $[m] := (m,m-1,\ldots,1)$ be
the corresponding shifted staircase shape.
Consider the {\em truncated} shifted staircase shape $[m]\setminus \kappa$,
where a partition $\kappa=(\kappa_1,...,\kappa_k)$, 
with $\kappa_i\leq m-i$ for all $1\le i\le k< m$, is deleted from the NE corner; 
namely, $\kappa_1$ cells are deleted from the (right) end of the first row, 
$\kappa_2$ cells are deleted from the end of the second row, etc. 
Let $N$ be the size of $[m]\setminus \kappa$. 
A standard Young tableau (SYT) of truncated shifted staircase shape 
$[m]\setminus \kappa$ is a labeling of the cells of this shape by $\{1,\ldots, N\}$ 
such that labels increase across rows (from left to right), down columns (from top
to bottom) and down the main diagonal (from top left to bottom right).


\medskip

\begin{example}\label{ex5} 
\rm \noindent There are four SYT of shape $[4]\setminus (1)$:
%

$$
{\vcenter
 {\offinterlineskip
 \halign{&\mystrut\vrule#&\mybox{\hss$\scriptstyle#$\hss}\cr
  \hr{7}\cr
  &1&\omit&2&\omit&3&& &\omit\cr
  \hr{3}&\nr{3}&\hr{3}&\nr{1}&\hr{1}\cr
  \omit& &&4&\omit&5&\omit&6&\cr
  \nr{2}&\hr{3}&\nr{3}&\hr{1}\cr
  \omit& &\omit& &&7&\omit&8&\cr
  \nr{4}&\hr{3}&\nr{1}&\hr{1}\cr
  \omit& &\omit& &\omit& &&9&\cr
  \nr{6}&\hr{3}\cr
   }}}
\ \ ,\
{\vcenter
 {\offinterlineskip
 \halign{&\mystrut\vrule#&\mybox{\hss$\scriptstyle#$\hss}\cr
  \hr{7}\cr
  &1&\omit&2&\omit&4&& &\omit\cr
  \hr{3}&\nr{3}&\hr{3}&\nr{1}&\hr{1}\cr
  \omit& &&3&\omit&5&\omit&6&\cr
  \nr{2}&\hr{3}&\nr{3}&\hr{1}\cr
  \omit& &\omit& &&7&\omit&8&\cr
  \nr{4}&\hr{3}&\nr{1}&\hr{1}\cr
  \omit& &\omit& &\omit& &&9&\cr
  \nr{6}&\hr{3}\cr
   }}}
\ \ ,\
{\vcenter
 {\offinterlineskip
 \halign{&\mystrut\vrule#&\mybox{\hss$\scriptstyle#$\hss}\cr
  \hr{7}\cr
  &1&\omit&2&\omit&3&& &\omit\cr
  \hr{3}&\nr{3}&\hr{3}&\nr{1}&\hr{1}\cr
  \omit& &&4&\omit&5&\omit&7&\cr
  \nr{2}&\hr{3}&\nr{3}&\hr{1}\cr
  \omit& &\omit& &&6&\omit&8&\cr
  \nr{4}&\hr{3}&\nr{1}&\hr{1}\cr
  \omit& &\omit& &\omit& &&9&\cr
  \nr{6}&\hr{3}\cr
   }}}
\ \ ,\
{\vcenter
 {\offinterlineskip
 \halign{&\mystrut\vrule#&\mybox{\hss$\scriptstyle#$\hss}\cr
  \hr{7}\cr
  &1&\omit&2&\omit&4&& &\omit\cr
  \hr{3}&\nr{3}&\hr{3}&\nr{1}&\hr{1}\cr
  \omit& &&3&\omit&5&\omit&7&\cr
  \nr{2}&\hr{3}&\nr{3}&\hr{1}\cr
  \omit& &\omit& &&6&\omit&8&\cr
  \nr{4}&\hr{3}&\nr{1}&\hr{1}\cr
  \omit& &\omit& &\omit& &&9&\cr
  \nr{6}&\hr{3}\cr
   }}}
\ \ .
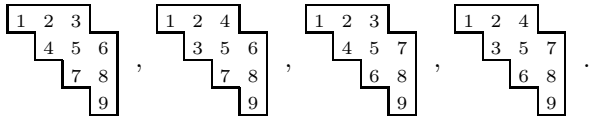
$$
\end{example}

\bigskip

Similarly, for nonnegative integers $m$ and $n$, let
$(n^m)=(n,\ldots,n)$ ($m$ parts) be the corresponding rectangular
shape. Consider the truncated rectangular shape
$(n^m)\setminus \kappa$, where $\kappa \subseteq (n^m)$
is deleted from the NE corner; namely,
$\kappa_1$ cells are deleted from the end of the first row,
$\kappa_2$ cells are deleted from the end of the second row, etc.
Letting $N$ be the size of $(n^m)\setminus \kappa$, a SYT of truncated shape
$(n^m)\setminus \kappa$ is a labeling of the cells of this shape
by $\{1,\ldots, N\}$, such that labels increase along rows
(from left to right) and down columns (from top to bottom).

\bigskip

Preliminary computer experiments hinted that a remarkable
phenomenon occurs when a {\em square} is truncated from a
staircase shape: while the number of SYT of truncated staircase
shape $[m]\setminus (k^k)$ increases exponentially as a function
of the size $N = {{m+1} \choose 2}-k^2$ of the shape, all the prime
factors of this number are actually smaller than the size. A
similar phenomenon occurs for squares truncated from rectangular
shapes \footnote{It should be noted that this factorization
phenomenon does not hold in the general case. For example, the
number of SYT of truncated rectangular shape $(7^6)\setminus (2)$,
of size $40$, has $5333$ as a prime factor.}. In this paper,
product formulas for these (and related) truncated shapes,
explaining the above factorization phenomenon, will be proved;
see Corollaries~\ref{t.stair_minus_sq_plus_1}, \ref{t.stair_minus_sq}, 
\ref{t.rect_minus_sq_plus_1} and \ref{t.rect_minus_sq} below.


\section{Pivoting}\label{sec:idea}

The basic idea of the proofs in the following sections is
to enumerate the SYT $T$ of a given shape $\zeta$ by mapping them bijectively
to pairs $(T_1,T_2)$ of SYT of some other shapes.
This will be done in two distinct (but superficially similar) ways, which will
complement each other and lead to the desired results.

The first bijection will be applied only to non-truncated shapes
of two types: a rectangle and a shifted staircase.
Let $N$ be the size of the shape $\zeta$ (of either of the aforementioned types), 
and fix an integer $1\le t\le N$.
Subdivide the entries in a SYT $T$ of shape $\zeta$ into
those that are less than or equal to $t$ and those that are greater than $t$.
The entries $\leq t$ constitute $T_1$. To obtain $T_2$, replace each entry
$i>t$ of $T$ by $N-i+1$, and suitably transpose (actually, reflect in the line $y=x$)
the resulting array. It is easy to see that both $T_1$ and $T_2$ are SYT.

This is illustrated schematically in the case of a shifted staircase shape
by the following shifted SYT, where $\zeta=(5,4,3,2,1)$, $N=15$ and $t=7$.
$$
T\ =\ {\vcenter
 {\offinterlineskip
 \halign{&\mystrut\vrule#&\mybox{\hss$\scriptstyle#$\hss}\cr
  \hr{11}\cr
  &1&\omit&2&\omit&3&\omit&6&&10&\cr
  \hr{3}&\nr{3}&\hr{3}&\nr{1}&\hr{1}\cr
  \omit& &&4&\omit&5&&8&\omit&11&\cr
  \nr{2}&\hr{3}&\nr{1}&\hr{1}&\nr{3}&\hr{1}\cr
  \omit& &\omit& &&7&&9&\omit&13&\cr
  \nr{4}&\hr{3}&\nr{3}&\hr{1}\cr
  \omit& &\omit& &\omit& &&12&\omit&14&\cr
  \nr{6}&\hr{3}&\nr{1}&\hr{1}\cr
  \omit& &\omit& &\omit& &\omit& &&15 &\cr
  \nr{8}&\hr{3}\cr
   }}}
\qquad\Leftrightarrow\qquad
T_1,T_2\ \  =\ \
{\vcenter
 {\offinterlineskip
 \halign{&\mystrut\vrule#&\mybox{\hss$\scriptstyle#$\hss}\cr
  \hr{11}\cr
  &1&\omit&2&\omit&3&\omit&6&\cr
  \hr{3}&\nr{3}&\hr{3}&\nr{1}&\hr{1}\cr
  \omit& &&4&\omit&5&\cr
  \nr{2}&\hr{3}&\nr{1}&\hr{1}&\nr{3}&\hr{1}\cr
  \omit& &\omit& &&7&\cr
  \nr{4}&\hr{3}&\nr{3}&\hr{1}\cr
}}}
\ ,\
{\vcenter
 {\offinterlineskip
 \halign{&\mystrut\vrule#&\mybox{\hss$\scriptstyle#$\hss}\cr
  \hr{11}\cr
  &1&\omit&2&\omit&3&\omit&5&\omit&6&\cr
  \hr{3}&\nr{5}&\hr{3}\cr
  \omit& &&4&\omit&7&\omit&8&\cr
  \nr{2}&\hr{7}\cr
   }}}
$$
In terms of shapes we have
$$
[\zeta]\ =\ {\vcenter
 {\offinterlineskip
 \halign{&\mystrut\vrule#&\mybox{\hss$\scriptstyle#$\hss}\cr
  \hr{11}\cr
  &\sigma&\omit& &\omit& &\omit& &&&\cr
  \hr{3}&\nr{3}&\hr{3}&\nr{1}&\hr{1}\cr
  \omit& && &\omit& && &\omit& &\cr
  \nr{2}&\hr{3}&\nr{1}&\hr{1}&\nr{3}&\hr{1}\cr
  \omit& &\omit& && && &\omit& &\cr
  \nr{4}&\hr{3}&\nr{3}&\hr{1}\cr
  \omit& &\omit& &\omit& && &\omit& &\cr
  \nr{6}&\hr{3}&\nr{1}&\hr{1}\cr
  \omit& &\omit& &\omit& &\omit& &&\tau'&\cr
  \nr{8}&\hr{3}\cr
   }}}
\qquad\Leftrightarrow\qquad
[\sigma],[\tau]\ \  =\ \
{\vcenter
 {\offinterlineskip
 \halign{&\mystrut\vrule#&\mybox{\hss$\scriptstyle#$\hss}\cr
  \hr{11}\cr
  &\sigma&\omit& &\omit& &\omit& &\cr
  \hr{3}&\nr{3}&\hr{3}&\nr{1}&\hr{1}\cr
  \omit& && &\omit& &\cr
  \nr{2}&\hr{3}&\nr{1}&\hr{1}&\nr{3}&\hr{1}\cr
  \omit& &\omit& && &\cr
  \nr{4}&\hr{3}&\nr{3}&\hr{1}\cr
}}}
\ ,\
{\vcenter
 {\offinterlineskip
 \halign{&\mystrut\vrule#&\mybox{\hss$\scriptstyle#$\hss}\cr
  \hr{11}\cr
  &\tau&\omit& &\omit& &\omit& &\omit& &\cr
  \hr{3}&\nr{5}&\hr{3}\cr
  \omit& && &\omit& &\omit& &\cr
  \nr{2}&\hr{7}\cr
   }}}
$$
where $\tau'$ denotes the conjugate of the (strict) partition $\tau$. It
follows immediately from this argument that
$$
g^\zeta = \sum_{{\sigma\subseteq\zeta} \atop {|\sigma|=t}}\ g^\sigma\ g^{(\zeta/\sigma)'}\,,
$$
where $\tau = (\zeta/\sigma)'$ is the conjugate of the shifted skew shape $\zeta/\sigma$.
The shape $[\tau]$ is formed by deleting the cells of $[\sigma]$ from
those of $[\zeta]$ and then reflecting in the line $y = x$.

In the case of a rectangular shape we have
$$
[\zeta]\ =\ {\vcenter
 {\offinterlineskip
 \halign{&\mystrut\vrule#&\mybox{\hss$\scriptstyle#$\hss}\cr
  \hr{11}\cr
  &\sigma &\omit& &\omit& &\omit& &\omit& &\cr
  \hr{1}&\nr{3}&\hr{7}\cr
  & &\omit&  && &\omit& &\omit& &\cr
  \hr{1}&\nr{1}&\hr{3}&\nr{5}&\hr{1}\cr
  & && &\omit& &\omit& &\omit&\tau'&\cr
  \hr{11}\cr
   }}}
\qquad\Leftrightarrow\qquad
[\sigma],[\tau]\ \  =\ \
{\vcenter
  {\offinterlineskip
 \halign{&\mystrut\vrule#&\mybox{\hss$\scriptstyle#$\hss}\cr
  \hr{11}\cr
  &\sigma &\omit& &\omit& &\omit& &\omit& &\cr
  \hr{1}&\nr{3}&\hr{7}\cr
  & &\omit&  &\cr
  \hr{1}&\nr{1}&\hr{3}\cr
  & &\cr
  \hr{3}\cr
   }}}
\ ,\
{\vcenter
 {\offinterlineskip
 \halign{&\mystrut\vrule#&\mybox{\hss$\scriptstyle#$\hss}\cr
  \hr{5}\cr
  &\tau &\omit& &\cr
  \hr{1}&\nr{3}&\hr{1}\cr
  & &\omit& &\cr
  \hr{1}&\nr{3}&\hr{1}\cr
  & &\omit& &\cr
  \hr{1}&\nr{1}&\hr{3}\cr
  & &\cr
  \hr{3}\cr
   }}}
$$
with
$$
f^\zeta = \sum_{\sigma\subseteq\zeta\atop |\sigma|=t}\ f^\sigma\, f^{(\zeta/\sigma)'}\,.
$$


\smallskip

The second bijection will be applied to truncated shapes,
where a partition is truncated from the NE corner of either a rectangle or a shifted staircase.
Given such a truncated shape $\zeta$, choose a {\em pivot cell} $P$. 
This is a cell of $\zeta$ which belongs to its NE boundary, namely such that
there is no cell of $\zeta$ which is both strictly north
and strictly east of the cell $P$.
If $T$ is a SYT of shape $\zeta$, let $t$ be the entry of $T$ in the pivot cell $P$.
Subdivide the other entries of $T$ into those that are (strictly) less than $t$ and those that are greater than $t$.
The entries less than $t$ constitute $T_1$. To obtain $T_2$,
replace each entry $i>t$ of $T$ by $N-i+1$, where $N$ is the
total number of entries in $T$, and suitably transpose the resulting array.
It is easy, once again, to see that both $T_1$ and $T_2$ are SYT.

By way of example,
consider the truncated rectangular shape specified by
$\zeta=(4,5,7,8,8)$ and let $P$ be the cell in position $(3,5)$.
Then for $t=17$ the map from a SYT $T$ of truncated shape
$[\zeta]$ to a corresponding pair $(T_1,T_2)$ is illustrated by
$$
T\ =\ {\vcenter
 {\offinterlineskip
 \halign{&\mystrut\vrule#&\mybox{\hss$\scriptstyle#$\hss}\cr
  \hr{9}\cr
  &1&\omit&2&\omit&4&\omit&9&\cr
  \hr{1}&\nr{7}&\hr{3}\cr
  &3&\omit&5&\omit&11&\omit&12&\omit&13&\cr
  \hr{1}&\nr{7}&\hr{7}\cr
  &6&\omit&8&\omit&14&\omit&15&&17&&21&\omit&24&\cr
  \hr{1}&\nr{3}&\hr{7}&\nr{3}&\hr{3}\cr
  &7&\omit&16&&18&\omit&20&\omit&25&\omit&26&\omit&27&\omit&30&\cr
  \hr{1}&\nr{1}&\hr{3}&\nr{11}&\hr{1}\cr
  &10&&19&\omit&22&\omit&23&\omit&28&\omit&29&\omit&31&\omit&32&\cr
  \hr{17}\cr
  }}}
\ \ \Leftrightarrow\ \
T_1,T_2\ =\
{\vcenter
 {\offinterlineskip
 \halign{&\mystrut\vrule#&\mybox{\hss$\scriptstyle#$\hss}\cr
  \hr{9}\cr
  &1&\omit&2&\omit&4&\omit&9&\cr
  \hr{1}&\nr{7}&\hr{3}\cr
  &3&\omit&5&\omit&11&\omit&12&\omit&13 &\cr
  \hr{1}&\nr{7}&\hr{3}\cr
  &6&\omit&8&\omit&14&\omit&15&\cr
  \hr{1}&\nr{3}&\hr{5}\cr
  &7&\omit&16&\cr
  \hr{1}&\nr{1}&\hr{3}\cr
  &10&\cr
  \hr{3}\cr
  }}}
\ \ ,\ \
 {\vcenter
 {\offinterlineskip
 \halign{&\mystrut\vrule#&\mybox{\hss$\scriptstyle#$\hss}\cr
  \hr{5}\cr
  &1&\omit&3&\cr
  \hr{1}&\nr{3}&\hr{3}\cr
  &2&\omit&6&\omit&9&\cr
  \hr{1}&\nr{5}&\hr{1}\cr
  &4&\omit&7&\omit&12&\cr
  \hr{1}&\nr{3}&\hr{3}\cr
  &5&\omit&8&\cr
  \hr{1}&\nr{3}&\hr{1}\cr
  &10&\omit&13&\cr
  \hr{1}&\nr{3}&\hr{1}\cr
  &11&\omit&15&\cr
  \hr{1}&\nr{1}&\hr{3}\cr
  &14&\cr
  \hr{3}\cr
    }}}
$$
In terms of shapes we have
$$
[\zeta]\ =\ {\vcenter
 {\offinterlineskip
 \halign{&\mystrut\vrule#&\mybox{\hss$\scriptstyle#$\hss}\cr
  \hr{9}\cr
  &\sigma &\omit& &\omit& &\omit& &\cr
  \hr{1}&\nr{7}&\hr{3}\cr
  & &\omit& &\omit& &\omit& &\omit& &\cr
  \hr{1}&\nr{7}&\hr{7}\cr
  & &\omit& &\omit& &\omit& &&P&& &\omit& &\cr
  \hr{1}&\nr{3}&\hr{7}&\nr{3}&\hr{3}\cr
  & &\omit& && &\omit& &\omit& &\omit& &\omit& &\omit& &\cr
  \hr{1}&\nr{1}&\hr{3}&\nr{11}&\hr{1}\cr
  & && &\omit& &\omit& &\omit& &\omit& &\omit& &\omit&\tau' &\cr
  \hr{17}\cr
  }}}
\ \ \Leftrightarrow\ \
[\sigma],[\tau]\ =\
{\vcenter
 {\offinterlineskip
 \halign{&\mystrut\vrule#&\mybox{\hss$\scriptstyle#$\hss}\cr
  \hr{9}\cr
  &\sigma &\omit& &\omit& &\omit& &\cr
  \hr{1}&\nr{7}&\hr{3}\cr
  & &\omit& &\omit& &\omit& &\omit& &\cr
  \hr{1}&\nr{7}&\hr{3}\cr
  & &\omit& &\omit& &\omit& &\cr
  \hr{1}&\nr{3}&\hr{5}\cr
  & &\omit& &\cr
  \hr{1}&\nr{1}&\hr{3}\cr
  & &\cr
  \hr{3}\cr
  }}}
\ \ ,\ \
 {\vcenter
 {\offinterlineskip
 \halign{&\mystrut\vrule#&\mybox{\hss$\scriptstyle#$\hss}\cr
  \hr{5}\cr
  &\tau &\omit& &\cr
  \hr{1}&\nr{3}&\hr{3}\cr
  & &\omit& &\omit& &\cr
  \hr{1}&\nr{5}&\hr{1}\cr
  & &\omit& &\omit& &\cr
  \hr{1}&\nr{3}&\hr{3}\cr
  & &\omit& &\cr
  \hr{1}&\nr{3}&\hr{1}\cr
  & &\omit& &\cr
  \hr{1}&\nr{3}&\hr{1}\cr
  & &\omit& &\cr
  \hr{1}&\nr{1}&\hr{3}\cr
  & &\cr
  \hr{3}\cr
    }}}
$$
where $\tau'$ is the conjugate of $\tau$.
%
Similarly for a truncated shifted staircase shape.

The particular shapes required in the sequel are the following:
\begin{equation}
\label{Eq-truncstairP}
{\vcenter
 {\offinterlineskip
 \halign{&\mystrut\vrule#&\mybox{\hss$\scriptstyle#$\hss}\cr
  \hr{17}\cr
  & &\omit& &\omit& &\omit& &\omit& &\omit& &\omit& &\omit& &\cr
  \hr{3}&\nr{13}&\hr{1}\cr
  \omit& && &\omit& &\omit& &\omit&\mu &\omit& &\omit& &\omit& &\cr
  \nr{2}&\hr{3}&\nr{11}&\hr{3}\cr
  \omit& &\omit& && &\omit& &\omit& &\omit& &\omit&  &\omit& &&P&\cr
  \nr{4}&\hr{19}\cr
  \omit& &\omit& &\omit& && &\omit&\lambda &\omit& &\omit& && && &\omit& &\omit& &\cr
  \nr{6}&\hr{3}&\nr{3}&\hr{3}&\nr{1}&\hr{1}&\nr{5}&\hr{1}\cr
  \omit& &\omit& &\omit& &\omit& && &\omit& && &\omit& && &\omit& &\omit& &\cr
  \nr{8}&\hr{3}&\nr{1}&\hr{1}&\nr{3}&\hr{1}&\nr{5}&\hr{1}\cr
  \omit& &\omit& &\omit& &\omit& &\omit& && && &\omit&\nu' && &\omit&\mu'  &\omit& &\cr
  \nr{10}&\hr{3}&\nr{3}&\hr{1}&\nr{5}&\hr{1}\cr
  \omit& &\omit& &\omit& &\omit& &\omit& &\omit& && &\omit& && &\omit& &\omit& &\cr
  \nr{12}&\hr{3}&\nr{1}&\hr{1}&\nr{5}&\hr{1}\cr
  \omit& &\omit& &\omit& &\omit& &\omit& &\omit& &\omit& && && &\omit& &\omit& &\cr
  \nr{14}&\hr{3}&\nr{5}&\hr{1}\cr
  \omit& &\omit& &\omit& &\omit& &\omit& &\omit& &\omit& &\omit& && &\omit& &\omit& &\cr
  \nr{16}&\hr{3}&\nr{3}&\hr{1}\cr
  \omit& &\omit& &\omit& &\omit& &\omit& &\omit& &\omit& &\omit& &\omit& && &\omit& &\cr
  \nr{18}&\hr{3}&\nr{1}&\hr{1}\cr
  \omit& &\omit& &\omit& &\omit& &\omit& &\omit& &\omit& &\omit& &\omit& &\omit& && &\cr
  \nr{20}&\hr{3}&\cr
}}}
\qquad\qquad\qquad
{\vcenter
 {\offinterlineskip
 \halign{&\mystrut\vrule#&\mybox{\hss$\scriptstyle#$\hss}\cr
  \hr{19}\cr
  & &\omit& &\omit& &\omit& &\omit& &\omit& &\omit& &\omit& &\omit& &\cr
  \hr{3}&\nr{15}&\hr{1}\cr
  \omit& && &\omit& &\omit& &\omit&\mu &\omit& &\omit& &\omit& &\omit& &\cr
  \nr{2}&\hr{3}&\nr{11}&\hr{7}\cr
  \omit& &\omit& && &\omit& &\omit& &\omit& &\omit&  &\omit& &&P&& &\omit& &\cr
  \nr{4}&\hr{15}&\nr{3}&\hr{1}\cr
  \omit& &\omit& &\omit& && &\omit&\lambda &\omit& &\omit& && && &\omit& &\omit& &\cr
  \nr{6}&\hr{3}&\nr{3}&\hr{3}&\nr{1}&\hr{1}&\nr{5}&\hr{1}\cr
  \omit& &\omit& &\omit& &\omit& && &\omit& && &\omit& && &\omit& &\omit& &\cr
  \nr{8}&\hr{3}&\nr{1}&\hr{1}&\nr{3}&\hr{1}&\nr{5}&\hr{1}\cr
  \omit& &\omit& &\omit& &\omit& &\omit& && && &\omit&\nu' && &\omit&\mu'  &\omit& &\cr
  \nr{10}&\hr{3}&\nr{3}&\hr{1}&\nr{5}&\hr{1}\cr
  \omit& &\omit& &\omit& &\omit& &\omit& &\omit& && &\omit& && &\omit& &\omit& &\cr
  \nr{12}&\hr{3}&\nr{1}&\hr{1}&\nr{5}&\hr{1}\cr
  \omit& &\omit& &\omit& &\omit& &\omit& &\omit& &\omit& && && &\omit& &\omit& &\cr
  \nr{14}&\hr{3}&\nr{5}&\hr{1}\cr
  \omit& &\omit& &\omit& &\omit& &\omit& &\omit& &\omit& &\omit& && &\omit& &\omit& &\cr
  \nr{16}&\hr{3}&\nr{3}&\hr{1}\cr
  \omit& &\omit& &\omit& &\omit& &\omit& &\omit& &\omit& &\omit& &\omit& && &\omit& &\cr
  \nr{18}&\hr{3}&\nr{1}&\hr{1}\cr
  \omit& &\omit& &\omit& &\omit& &\omit& &\omit& &\omit& &\omit& &\omit& &\omit& && &\cr
  \nr{20}&\hr{3}&\cr
  }}}
\end{equation}
and
\begin{equation}
\label{Eq-truncrectP}
{\vcenter
 {\offinterlineskip
 \halign{&\mystrut\vrule#&\mybox{\hss$\scriptstyle#$\hss}\cr
  \hr{13}\cr
  & &\omit& &\omit& &\omit& &\omit& &\omit& &\cr
  \hr{1}&\nr{11}&\hr{1}\cr
  & &\omit&\mu+\alpha &\omit& &\omit& &\omit& &\omit& &\cr
  \hr{1}&\nr{11}&\hr{3}\cr
  & &\omit& &\omit& &\omit& &\omit& &\omit& &&P&\cr
  \hr{19}\cr
  & &\omit&\lambda &\omit& &\omit& && &\omit& && &\omit& &\omit& &\cr
  \hr{1}&\nr{5}&\hr{3}&\nr{3}&\hr{1}&\nr{5}&\hr{1}\cr
  & &\omit& &\omit& && &\omit&\nu' &\omit& && &\omit&(\mu+\beta)' &\omit& &\cr
  \hr{7}&\nr{5}&\hr{1}&\nr{5}&\hr{1}\cr
    & &\omit& &\omit& &\omit& &\omit& &\omit& && &\omit& &\omit& &\cr
  \hr{19}\cr
}}}
\qquad\qquad\qquad\qquad
{\vcenter
 {\offinterlineskip
 \halign{&\mystrut\vrule#&\mybox{\hss$\scriptstyle#$\hss}\cr
  \hr{15}\cr
  & &\omit& &\omit& &\omit& &\omit& &\omit& &\omit& &\cr
  \hr{1}&\nr{13}&\hr{1}\cr
  & &\omit&\mu+\alpha &\omit& &\omit& &\omit& &\omit& &\omit& &\cr
  \hr{1}&\nr{11}&\hr{7}\cr
  & &\omit& &\omit& &\omit& &\omit& &\omit& && P && &\omit& &\cr
  \hr{15}&\nr{3}&\hr{1}\cr
  & &\omit&\lambda &\omit& &\omit& && &\omit& && &\omit& &\omit& &\cr
  \hr{1}&\nr{5}&\hr{3}&\nr{3}&\hr{1}&\nr{5}&\hr{1}\cr
  & &\omit& &\omit& && &\omit&\nu' &\omit& && &\omit& (\mu+\beta)' &\omit& &\cr
  \hr{7}&\nr{5}&\hr{1}&\nr{5}&\hr{1}\cr
    & &\omit& &\omit& &\omit& &\omit& &\omit& && &\omit& &\omit& &\cr
  \hr{19}\cr
}}}
\end{equation}
A crucial property of these particular shapes is that their subdivision
gives rise to shapes $\mu \cup \lambda$ and $\mu \cup \nu$
(or $(\mu + \alpha) \cup \lambda$ and $(\mu + \beta) \cup \nu$)
which are {\em not} truncated.

A more explicit relation between $\lambda$ and $\nu$ will be given later.
For the time being it suffices to observe that
\begin{equation}
\label{Eq-gzetaP}
g^\zeta = \sum_{\lambda\subseteq[m]}\ g^{\mu\cup\lambda}\ g^{\mu\cup([m]/\lambda)'}
\end{equation}
and
\begin{equation}
\label{Eq-fzetaP}
f^\zeta = \sum_{\lambda\subseteq(n^m)}\ f^{(\mu+\alpha)\cup\lambda}\ f^{(\mu+\beta)\cup((m^n)/\lambda)'}\,.
\end{equation}
Since $\zeta$ is a truncated shape, the notation on the LHS of the above equalities
is to be taken in a generalized sense.

\section{Truncated Shifted Staircase Shapes}\label{sec:staircase}

In this section, $\la = (\la_1, \ldots, \la_m)$ (with $\la_1 > \ldots > \la_m > 0$ integers)
will be a strict partition,
with $g^{\la}$ denoting the number of SYT of shifted shape $\la$.
We shall use the union operation on partitions, defined in Section~\ref{sec:basic}.

For any nonnegative integer $m$, let $[m] := (m,m-1,\ldots,1)$ be
the corresponding shifted staircase shape.  Schur's formula
(Proposition ~\ref{t.num_shifted_SYT}) implies the following.

\begin{observation}\label{t.shifted_staircase}
The number of SYT of shifted staircase shape $[m]$ is
$$
g^{[m]} = M! \cdot \prod_{i=0}^{m-1} \frac{i!}{(2i+1)!},
$$
where $M := |[m]| = {m+1 \choose 2}$.
\end{observation}

\medskip

We shall now use the first bijection outlined in Section~\ref{sec:idea}.

\begin{lemma}\label{t.complement_shifted}
Let $m$ and $t$ be nonnegative integers, with $t \le {m+1 \choose 2}$.
Let $T$ be a SYT of shifted staircase shape $[m]$,
let $T_1$ be the set of all cells in $T$ with values at most $t$,
and let $T_2$ be obtained from $T \setminus T_1$ by transposing the shape and
replacing each entry $i$ by $M-i+1$. Then:
\begin{itemize}
\item[(i)]
$T_1$ and $T_2$ are shifted SYT.
\item[(ii)]
Treating strict partitions as sets, $[m]$ is the disjoint union of
the shape of $T_1$ and the the shape of $T_2$.
\end{itemize}
\end{lemma}

\noindent{\bf Proof.} (i) is clear. In order to prove (ii), denote
the shifted shapes of $T_1$ and $T_2$ by $\la^1$ and $\la^2$,
respectively. The borderline between $T_1$ and $T \setminus T_1$ is a
lattice path of length exactly $m$, starting at the NE corner of
the staircase shape $[m]$
and using only S and W steps. 
If the first step is S then the first part of $\la^1$ is $m$, and
the rest (of both $\la^1$ and $\la^2$) corresponds to a lattice
path in $[m-1]$. Similarly, if the first step is W then the first
part of $\la^2$ is $m$, and the rest corresponds to a lattice path
in $[m-1]$. Thus exactly one of $\la^1$, $\la^2$ has a part equal
to $m$, and the whole result follows by induction on $m$.

\qed


\begin{corollary}\label{t.sum_shifted}
For any nonnegative integers $m$ and $t$ with $t \le {m+1 \choose 2}$,
$$
\sum_{\la \subseteq [m]\atop |\la|=t} g^{\la}g^{\la^c} = g^{[m]}.
$$
Here summation is over all strict partitions $\la$ with the prescribed restrictions,
and $\la^c$ is the complement of $\la$ in $[m]$
(where strict partitions are treated as sets).
In particular, the LHS is independent of $t$.
\end{corollary}

\begin{lemma}\label{t.coeff_shifted}
Let $\la$ and $\la^c$ be strict partitions whose disjoint union (as sets) is $[m]$,
and let $\mu=(\mu_1,\ldots,\mu_k)$ with $\mu_1 > \ldots > \mu_k >m$.
Let $\hg^{\la} := g^{\la}/|\la|!$.
Then
$$
\hg^{\mu \cup \la} \hg^{\mu \cup \la^c} \hg^{[m]} =
\hg^{\la} \hg^{\la^c} \hg^{\mu \cup [m]} \hg^{\mu}.
$$
Equivalently,
$$
g^{\mu \cup \la} g^{\mu \cup \la^c} = c(\mu,|\la|,|\la^c|) \cdot g^{\la} g^{\la^c},
$$
where
$$
c(\mu,|\la|,|\la^c|) = 
\frac{g^{\mu \cup [m]} g^{\mu}}{g^{[m]}} \cdot 
\frac{|[m]|!(|\mu|+|\la|)!(|\mu|+|\la^c|)!}{(|\mu|+|[m]|)!|\mu|!|\la|!|\la^c|!}
$$
depends only on the sizes $|\la|$ and $|\la^c|$ and not on the actual partitions 
$\la$ and $\la^c$.
\end{lemma}
\noindent{\bf Proof.}
By Proposition~\ref{t.num_shifted_SYT},
\begin{equation}\label{e.hg_identity}
\frac{\hg^{\mu \cup \la} \hg^{\mu \cup \la^c}}{\hg^{\la} \hg^{\la^c}} =
\left(\prod_i \frac{1}{\mu_i!}\prod_{i<j} \frac{\mu_i-\mu_j}{\mu_i+\mu_j}\right)^2
\cdot 
\prod_{i,j} \frac{\mu_i-\la_j}{\mu_i+\la_j}
\prod_{i,j} \frac{\mu_i-\la^c_j}{\mu_i+\la^c_j}.
\end{equation}
By the assumption on $\la$ and $\la^c$,
$$
\prod_{j} \frac{\mu_i-\la_j}{\mu_i+\la_j}
\prod_{j} \frac{\mu_i-\la^c_j}{\mu_i+\la^c_j}
= \prod_{j=1}^{m} \frac{\mu_i-j}{\mu_i+j}\qquad(\forall i).
$$
Thus the RHS of~(\ref{e.hg_identity}) is independent of $\la$ and $\la^c$.
Substituting $\la = [m]$ (and $\la^c = [0]$, the empty partition) yields
$$
\frac{\hg^{\mu \cup \la} \hg^{\mu \cup \la^c}}{\hg^{\la} \hg^{\la^c}} =
\frac{\hg^{\mu \cup [m]} \hg^{\mu}}{\hg^{[m]}},
$$
which is the desired identity. The other equivalent formulation follows readily.

\qed

\bigskip


A technical lemma, which will be used to prove
Theorems~\ref{main1111} and~\ref{main2222},
is the following.

\begin{lemma}\label{t.binomial}
Let $t_1$, $t_2$ and $N$ be nonnegative integers. Then
$$
\sum_{i=0}^{N} {t_1 + i \choose t_1} {t_2 + N - i \choose t_2} =
{t_1 + t_2 + N + 1 \choose t_1 + t_2 + 1}.
$$
\end{lemma}

\noindent{\bf Proof.} This is a classical binomial identity, which
follows for example from computation of the coefficients of $x^N$
on both sides of the identity
$$
(1-x)^{-(1+t_1)} \cdot (1-x)^{-(1+t_2)} = (1-x)^{-(2+t_1+t_2)}.
$$

\qed

\bigskip

\begin{theorem}\label{main1111}
Let $m$ be a nonnegative integer, denote $M := {m+1 \choose 2}$,
and let $\mu=(\mu_1,\ldots,\mu_k)$ be a strict partition with $\mu_1 > \ldots > \mu_k >m$.
Then
$$
\sum_{\la \subseteq [m]} g^{\mu \cup \la}g^{\mu \cup \la^c}
=
g^{\mu \cup [m]} g^{\mu} \cdot \frac{(M+2|\mu|+1)!|\mu|!}{(M+|\mu|)!(2|\mu|+1)!}.
$$
\end{theorem}
\noindent{\bf Proof.}
Restrict the summation on the LHS to strict partitions $\la$ of a fixed size $|\la|=t$.
By Lemma~\ref{t.coeff_shifted} and Corollary~\ref{t.sum_shifted},
$$
\sum_{\la \subseteq [m] \atop |\la|=t} g^{\mu \cup \la}g^{\mu \cup \la^c} =
c(\mu,t,M-t) \cdot \sum_{\la \subseteq [m] \atop |\la|=t} g^{\la}g^{\la^c} =
c(\mu,t,M-t) \cdot g^{[m]}.
$$
Now sum over all $t$ and use the explicit formula for $c(\mu,t,M-t)$ (from
Lemma~\ref{t.coeff_shifted}) together with Lemma~\ref{t.binomial}: 
\begin{eqnarray*}
\sum_{\la \subseteq [m]} g^{\mu \cup \la}g^{\mu \cup \la^c} 
&=&
g^{\mu \cup [m]} g^{\mu} \cdot \frac{M!}{(|\mu|+M)!|\mu|!} \cdot
\sum_{t=0}^{M} \frac{(|\mu|+t)!(|\mu|+M-t)!}{t!(M-t)!}\\
&=&
g^{\mu \cup [m]} g^{\mu} \cdot {{|\mu|+M} \choose {|\mu|}}^{-1} \cdot
\sum_{t=0}^{M} {{|\mu|+t} \choose {|\mu|}} {{|\mu|+M-t} \choose {|\mu|}}\\
&=&
g^{\mu \cup [m]} g^{\mu} \cdot {{|\mu|+M} \choose {|\mu|}}^{-1} \cdot
{{2|\mu|+M+1} \choose {2|\mu|+1}}.
\end{eqnarray*}

\qed

\medskip

We shall apply this theorem to several special cases.
In each case the result will follow from an application of equation~(\ref{Eq-gzetaP})
to one or the other of the shapes illustrated in diagram~(\ref{Eq-truncstairP}),
where the union of $\la$ and $\nu$ is the shifted staircase shape $[m]$.
Using Lemma~\ref{t.complement_shifted},
we can now state explicitly the relation between these partitions, mentioned before
equation~(\ref{Eq-gzetaP}): 
$\nu=([m]/\lambda)'=\lambda^c$.

\smallskip

First, take $\mu=(m+k,\ldots,m+1)$ ($k$ parts), for $k\ge 1$. 
This corresponds to truncating a $k \times k$ square from the NE corner of 
a shifted staircase shape $[m+2k]$, but adding back the SW corner of this square; 
see the first shape in diagram (\ref{Eq-truncstairP}).

\begin{corollary}\label{t.stair_minus_sq_plus_1}
The number of SYT of truncated shifted staircase shape $[m+2k] \setminus (k^{k-1},k-1)$ is
$$
g^{[m+k]} g^{(m+k,\ldots,m+1)} \cdot \frac{N!|\mu|!}{(N-|\mu|-1)!(2|\mu|+1)!},
$$
where $N = {m+2k+1 \choose 2} - k^2 +1$ is the size of the shape and $|\mu| = k(2m+k+1)/2$.
\end{corollary}

The special case $k=1$ (with $\mu=(m+1)$) gives back
the number $g^{[m+2]}$ of SYT of shifted staircase shape $[m+2]$:
$$
g^{[m+1]} \cdot \frac{N!(m+1)!}{(N-m-2)!(2m+3)!} 
= 
N! \cdot \prod_{i=0}^{m+1} \frac{i!}{(2i+1)!}
=
g^{[m+2]},
$$
where $N = (m+2)(m+3)/2$ is the size of the shape.
This agrees, of course, with Observation~\ref{t.shifted_staircase}.

The special case $k=2$ (with $\mu=(m+2,m+1)$) corresponds to
truncating a small shifted staircase shape $[2] = (2,1)$ from the
shifted staircase shape $[m+4]$. Thus,
the number of SYT of truncated shifted staircase shape $[m+4] \setminus [2]$ is
\begin{eqnarray*}
& & 
g^{[m+2]} g^{(m+2,m+1)} \cdot \frac{N!(2m+3)!}{(N-2m-4)!(4m+7)!}\\
&=&
\prod_{i=0}^{m+1} \frac{i!}{(2i+1)!} \cdot \frac{(2m+3)!}{(m+2)!(m+1)!(2m+3)} \cdot \frac{N!(2m+3)!}{(4m+7)!} \\
&=&
N! \cdot \frac{2}{(4m+7)!(m+2)} \cdot \prod_{i=0}^{m-1} \frac{i!}{(2i+1)!},
\end{eqnarray*}
where $N = {{m+5} \choose 2} - 3 = (m+2)(m+7)/2$ is the size of the shape.

\medskip

Now take $\mu=(m+k+1,\ldots,m+3,m+1)$ ($k$ parts), for $k\ge 2$. This corresponds to
truncating a $(k-1) \times (k-1)$ square from the NE corner of a shifted
staircase shape $[m+2k]$;
see the second shape in diagram (\ref{Eq-truncstairP}).

\begin{corollary}\label{t.stair_minus_sq}
The number of SYT of truncated shifted staircase shape $[m+2k] \setminus ((k-1)^{k-1})$ is
$$
g^{(m+k+1,\ldots,m+3,m+1,\ldots,1)} g^{(m+k+1,\ldots,m+3,m+1)} \cdot \frac{N!|\mu|!}{(N-|\mu|-1)!(2|\mu|+1)!},
$$
where $N = {m+2k+1 \choose 2} - (k-1)^2$ is the size of the shape and $|\mu| = k(2m+k+3)/2 - 1$.
\end{corollary}

In particular, take $k=2$ and $\mu=(m+3,m+1)$. 
This corresponds to truncating the NE corner cell of a shifted staircase shape $[m+4]$.
The corresponding enumeration problem was actually the original motivation for the current work, 
because of its combinatorial interpretation, as explained in~\cite{AR_tft2}.
Thus, the number of SYT of truncated shifted staircase shape 
$[m+4]\setminus (1)$ is
\begin{eqnarray*}
& &
g^{(m+3,m+1,\ldots,1)} g^{(m+3,m+1)} \cdot \frac{N!(2m+4)!}{(N-2m-5)!(4m+9)!}\\
&=&
\frac{N! \cdot 4(2m+3)}{(4m+9)! \cdot (m+3)} \cdot \prod_{i=0}^{m-1} \frac{i!}{(2i+1)!},
\end{eqnarray*}
where $N = {{m+5} \choose 2} - 1 = (m+3)(m+6)/2$ is the size of the shape.
For $m = 0$ ($N = 9$) this number is $4$, as shown in Example~\ref{ex5}.




\bigskip

\section{Truncated Rectangular Shapes}\label{sec:rectangular}

In this section, $\la = (\la_1, \ldots, \la_m)$ (with $\la_1 \ge \ldots \ge \la_m \ge 0$ integers)
will be a partition with (at most) $m$ parts.
Two partitions which differ only in trailing zeros will be considered equal.
Denote by $f^{\la}$ the number of SYT of regular (non-shifted) shape $\la$.

For any nonnegative integers $m$ and $n$, let $(n^m) :=
(n,\ldots,n)$ ($m$ times) be the corresponding rectangular shape.
The Frobenius-Young formula (Proposition~\ref{t.num_regular_SYT})
implies the following.

\begin{observation}\label{t.rectangle}
The number of SYT of rectangular shape $(n^m)$ is
$$
f^{(n^m)} = (mn)! \cdot \frac{F_m F_n}{F_{m+n}},
$$
where
$$
F_m := \prod_{i=0}^{m-1} i!.
$$
\end{observation}

Recall from Section~\ref{sec:basic} the definition of 
the sum $\la+\mu$ of two partitions $\la$ and $\mu$. 
Note that if either $\la$ or $\mu$ is a strict partition then so is $\la+\mu$.

\begin{lemma}\label{t.complement_regular}
Let $m$, $n$ and $t$ be nonnegative integers, with $t \le mn$.
Let $T$ be a SYT of rectangular shape $(n^m)$,
let $T_1$ be the set of all cells in $T$ with values at most $t$,
and let $T_2$ be obtained from $T \setminus T_1$ by transposing the shape and
replacing each entry $i$ by $mn-i+1$. Then:
\begin{itemize}
\item[(i)]
$T_1$ and $T_2$ are SYT. 
\item[(ii)]
Denote by $\la^1$ and $\la^2$
the shapes of $T_1$ and $T_2$, respectively, and treat strict
partitions as sets. Then the strict partition $[m+n]$ is the
disjoint union of the strict partitions $\la^1 + [m]$ and $\la^2 + [n]$.
\end{itemize}
\end{lemma}

\noindent{\bf Proof.} 
(i) is clear; let us prove (ii). 
The borderline between $T_1$ and $T \setminus T_1$ is a lattice path of length
exactly $m+n$, starting at the NE corner of the rectangular shape $(n^m)$,
using only S and W steps, and ending at the SW corner of this shape. 
If the first step is S then the first part of $\la^1 + [m]$ is
$m+n$, and the rest (of both $\la^1 + [m]$ and $\la^2 + [n]$)
corresponds to a lattice path in $n^{m-1}$. Similarly, if the
first step is W then the first part of $\la^2 + [n]$ is $m+n$, and
the rest corresponds to a lattice path in $(n-1)^m$. Thus exactly
one of the strict partitions $\la^1 + [m]$ and $\la^2 + [n]$ has a
part equal to $m+n$, and the whole result follows by induction on
$m+n$.

\qed


\begin{corollary}\label{t.sum_regular}
For any nonnegative integers $m$, $n$ and $t$ with $t \le mn$,
$$
\sum_{\la \subseteq (n^m)\atop |\la|=t} f^{\la}f^{\la^c} = f^{(n^m)}.
$$
Here summation is over all partitions $\la$ with the prescribed restrictions,
and $\la^c$ is such that $\la^c + [n]$ is the complement of $\la + [m]$ in $[m+n]$
(where strict partitions are treated as sets).
In particular, the LHS is independent of $t$.
\end{corollary}

\begin{lemma}\label{t.coeff_regular}
Let $\la$ and $\la^c$ be partitions such that
$[m+n]$ is the disjoint union of $\la + [m]$ and $\la^c + [n]$,
and let $\mu=(\mu_1,\ldots,\mu_k)$ be an arbitrary partition
($\mu_1 \ge \ldots \ge \mu_k \ge 0$).
Let $\hf^{\la} := f^{\la}/|\la|!$.
Then
$$
\hf^{(\mu + (n^k)) \cup \la} \hf^{(\mu + (m^k)) \cup \la^c} = 
\hf^{\la} \hf^{\la^c} \hf^{\mu + (m+n)^k} \hf^{\mu}.
$$
Equivalently,
$$
f^{(\mu + (n^k)) \cup \la} f^{(\mu + (m^k)) \cup \la^c} = d(\mu,|\la|,|\la^c|) \cdot f^{\la} f^{\la^c},
$$
where
$$
d(\mu,|\la|,|\la^c|) = f^{\mu + ((m+n)^k)} f^{\mu} \cdot
                       \frac{(|\mu|+nk+|\la|)!(|\mu|+mk+|\la^c|)!}{(|\mu|+(m+n)k)!(|\mu|)!(|\la|)!(|\la^c|)!}.
$$
\end{lemma}
\noindent{\bf Proof.}
From the assumptions it follows that $\la$ is contained in $(n^m)$.
We may thus assume that it has $m$ (nonnegative) parts.
Similarly, $\la^c$ is contained in $(m^n)$ and we may assume that it has $n$ (nonnegative) parts.
Thus $((\mu + (n^k)) \cup \la$ has $k+m$ parts and $(\mu + (m^k)) \cup \la^c$ has $k+n$ parts.
By Proposition~\ref{t.num_regular_SYT},
\begin{eqnarray*}
& &\frac{\hf^{(\mu + (n^k)) \cup \la} \hf^{(\mu + (m^k)) \cup \la^c}}{\hf^{\la} \hf^{\la^c}} = \\
& & \cdot \left( \prod_{i=1}^{k} \frac{1}{(\mu_i+m+n+k-i)!}\prod_{i<j} (\mu_i-\mu_j-i+j) \right)^2
\cdot\\
& & \cdot
\prod_{i,j} (\mu_i+n-\la_j-i+k+j)
\prod_{i,j} (\mu_i+m-\la^c_j-i+k+j).
\end{eqnarray*}
By the assumption on $\la$ and $\la^c$,
$$
\prod_{j} (\mu_i+n-\la_j-i+k+j)
\prod_{j} (\mu_i+m-\la^c_j-i+k+j)
= \prod_{j=1}^{m+n} (\mu_i-i+k+j)\qquad(\forall i).
$$
Since
$$
\frac{1}{(\mu_i+m+n+k-i)!} \cdot \prod_{j=1}^{m+n} (\mu_i-i+k+j)
= \frac{1}{(\mu_i+k-i)!}\qquad(\forall i),
$$
an application of Proposition~\ref{t.num_regular_SYT} to $f^{\mu}$ and to $f^{\mu+((m+n)^k)}$
gives the desired result.

\qed

\begin{theorem}\label{main2222}
Let $m$ and $n$ be nonnegative integers,
and let $\mu=(\mu_1,\ldots,\mu_k)$ be a partition ($\mu_1 \ge \ldots \ge \mu_k \ge 0$).
Then
\begin{eqnarray*}
& & \sum_{\la \subseteq (n^m)} f^{(\mu + (n^k)) \cup \la} f^{(\mu + (m^k)) \cup \la^c} =
f^{\mu + ((m+n)^k)} f^{\mu} f^{(n^m)} \cdot\\
& & \cdot {mn+2|\mu|+mk+nk+1 \choose mn} \cdot
\frac{(|\mu|+mk)!(|\mu|+nk)!}{(|\mu|+mk+nk)!(|\mu|)!}.
\end{eqnarray*}
\end{theorem}
\noindent{\bf Proof.}
Restrict the summation to partitions $\la$ of a fixed size $|\la|=t$.
By Lemma~\ref{t.coeff_regular} and Corollary~\ref{t.sum_regular},
$$
\sum_{\la \subseteq (n^m) \atop |\la|=t} f^{(\mu + (n^k)) \cup \la} f^{(\mu + (m^k)) \cup \la^c} =
d(\mu,t,M-t) \cdot \sum_{\la \subseteq (n^m) \atop |\la|=t} f^{\la} f^{\la^c} =
d(\mu,t,M-t) \cdot f^{(n^m)}.
$$
Now sum over all $t$ and use the explicit formula for $d(\mu,t,M-t)$ (from
Lemma~\ref{t.coeff_regular}) together with Lemma~\ref{t.binomial}, to obtain
the explicit formula above.

\qed

\medskip

Again, we shall apply this theorem in several special cases.
In each case the result is a special case of equation (\ref{Eq-fzetaP}) with
appropriately chosen $\zeta$ of one or the other of the shapes illustrated in
diagram~(\ref{Eq-truncrectP}). Note that, by Lemma~\ref{t.complement_regular},
$\nu=((n^m)/\lambda)'=\lambda^c$.

\smallskip

First, 
let $\alpha = (n^k)$, $\beta = (m^k)$ and $\mu = (0^k)$ (the empty partition with $k$ ``parts'').
This corresponds to truncating a $k \times k$ square from the NE corner 
of a rectangular shape $((n+k)^{m+k})$, but adding back the SW corner of this square;
see the first shape in diagram~(\ref{Eq-truncrectP}).

\begin{corollary}\label{t.rect_minus_sq_plus_1}
The number of SYT of truncated rectangular shape $((n+k)^{m+k}) \setminus (k^{k-1},k-1)$ is
$$
f^{((m+n)^k)} f^{(n^m)} \cdot {mn+mk+nk+1 \choose mn} \cdot \frac{(mk)!(nk)!}{(mk+nk)!}
=
\frac{N!(mk)!(nk)!}{(mk+nk+1)!} \cdot \frac{F_m F_n F_k}{F_{m+n+k}},
$$
where $N = (m+k)(n+k)-k^2+1 = mn+mk+nk+1$ is the size of the shape
and $F_n$ is as in Observation~\ref{t.rectangle}.
\end{corollary}

For $k=1$ we obtain
$$
f^{(n+1)^{m+1}} = \frac{N!m!n!}{(m+n+1)!} \cdot \frac{F_m F_n}{F_{m+n+1}} 
= N! \cdot \frac{F_{m+1} F_{n+1}}{F_{m+n+2}},
$$
in accordance with Observation~\ref{t.rectangle}.

For $k=2$ we obtain that
the number of SYT of truncated rectangular shape $((n+2)^{m+2}) \setminus (2,1)$ is
$$
\frac{N!(2m)!(2n)!}{(2m+2n+1)!} \cdot \frac{F_m F_n}{F_{m+n+2}},
$$
where $N = (m+2)(n+2)-3 = mn+2m+2n+1$ is the size of the shape.

\medskip

Now take $\alpha = (n^k)$, $\beta = (m^k)$ and $\mu = (1^{k-1},0)$, for $k\ge 2$.
This corresponds to truncating a $(k-1) \times (k-1)$ square from 
the NE corner of a rectangular shape $((n+k)^{m+k})$;
see the second shape in diagram~(\ref{Eq-truncrectP}).

\begin{corollary}\label{t.rect_minus_sq}
The number of SYT of truncated rectangular shape $((n+k)^{m+k}) \setminus ((k-1)^{k-1})$ is
$$
f^{((m+n+1)^{k-1},m+n)} f^{(n^m)} \cdot {mn+mk+nk+2k-1 \choose mn} \cdot
\frac{(mk+k-1)!(nk+k-1)!}{(mk+nk+k-1)!(k-1)!} =
$$
$$
=
\frac{N!(mk+k-1)!(nk+k-1)!(m+n+1)!k}{(mk+nk+2k-1)!} \cdot
\frac{F_m F_n F_{k-1}}{F_{m+n+k+1}},
$$
where $N = (m+k)(n+k)-(k-1)^2 = mn+mk+nk+2k-1$ is the size of the shape and $F_n$
is as in Observation~\ref{t.rectangle}.
\end{corollary}

In particular, letting $k=2$ and $\mu=(1,0)$ implies that
the number of SYT of truncated rectangular shape $((n+2)^{m+2})
\setminus (1)$ is
$$
\frac{N!(2m+1)!(2n+1)! \cdot 2}{(2m+2n+3)!(m+n+2)} \cdot
\frac{F_m F_n}{F_{m+n+2}},
$$
where $N = (m+2)(n+2)-1 = mn+2m+2n+3$ is the size of the shape and
$F_n$ is as in Observation~\ref{t.rectangle}.

\section{Final Remarks and Open Problems}\label{sec:final}

Most recently, the numbers of SYT of some other truncated shapes 
have been shown to have nice product formulas. 
Using complex analysis and volume computations, G.\ Panova proved 
such a formula for
a rectangular shape minus a staircase in the NE corner~\cite{Panova}. 
Computer experiments indicate that there are other shapes 
with a similar property.

\medskip

The number of SYT of a rectangular shape $(n^m)$ ($n\ne m$), 
with {\em two} cells truncated from the NE corner, seems in general to have large prime factors.
Nevertheless, the situation is different for $n=m$. 

\begin{conjecture} For $n\ge 2$
$$
f^{(n^n)\setminus (2)}=
\frac{ (n^2-2)! (3n-4)!^2 \cdot 6}{ (6n-8)! (2n-2)! (n-2)!^2} \cdot \frac{F_{n-2}^2}{F_{2n-4}},
$$
where
$$
F_m := \prod_{i=0}^{m-1} i!.
$$
\end{conjecture}

This empirical formula has been verified for $n\le 13$.

\smallskip

On the basis of further computations we conjecture that,
if either $ \kappa = (k^{k-1})$ or $\kappa = (k^{k-2},k-1)$ then,
for $n$ large enough,
both $f^{(n^n)\backslash\kappa}$ and $g^{[n]\backslash\kappa}$ 
have no prime factor larger than the size of the shape.

\begin{question}$\;$
\begin{itemize}
\item[(i)] 
Find and characterize the truncated (and other non ``unusual'') shapes
for which the number of SYT has no prime factor larger than the size of the shape.
\item[(ii)] 
Give explicit product formulas for the number of SYT of such shapes.
\end{itemize}
\end{question}

%
%
%
%
%

\bigskip

\noindent{\bf Acknowledgements:} 
Many thanks to Amitai Regev and Richard Stanley for stimulating comments and references. 
Thanks also to an anonymous referee for helpful comments. 
Extensive computer experiments have been conducted using the Maple and Sage
computer algebra systems.

This paper addresses a question presented by the authors in 
the Workshop on Combinatorics 
held at the Weizmann Institute in June 2010,
in honor of Doron Zeilberger's $60^{th}$ birthday. 
Doron, thanks for many years of inspiration!


\begin{thebibliography}{ABR4}

\bibitem{AR_tft2}
R.\ M.\ Adin and Y.\ Roichman,
{\em Triangle-free triangulations, hyperplane arrangements and shifted tableaux},
preprint 2010, arXiv:1009.2628.

\bibitem{A}
A.\ C.\ Aitken,
{\em The monomial expansion of determinantal symmetric functions},
Proc.\ Royal Soc.\ Edinburgh (A)~{\bf 61} (1943), 300--310.

\bibitem{FRT}
J.\ S.\ Frame, G.\ de B.\ Robinson and R.\ M.\ Thrall,
{\em The hook graphs of the symmetric group},
Canad.\ J.\ Math.~{\bf 6} (1954), 316--324.

\bibitem{F}
F.\ G.\ Frobenius,
{\em Uber die charaktere der symmetrischer gruppe},
Preuss.\ Akad.\ Wiss.\ Sitz.\ (1900), 516--534.

\bibitem{M}
P.\ A.\ MacMahon,
Combinatory Analysis, vols.\ 1 and 2,
Cambridge Univ.\ Press, London/New York, 1915, 1916;
reprinted in one volume by Chelsea, New York, 1960.

\bibitem{Md} I.\ G.\ Macdonald,
{\em Symmetric Functions and Hall Polynomials}, second ed., Oxford
Math.\ Monographs, Oxford Univ.\ Press, Oxford, 1995.

\bibitem{Panova} 
G.\ Panova, {\em Truncated tableaux and plane partitions},
preprint 2010, arXiv:1011.0795.

\bibitem{Schur} I.\ Schur, {\em On the representation of the symmetric and alternating groups by fractional linear
substitutions}, Translated from the German [J. Reine Angew. Math.
{\bf 139} (1911), 155–-250] by Marc-Felix Otto. Internat. J.
Theoret. Phys. {\bf 40} (2001), 413–-458.

\bibitem{St_EC2}
R.\ P.\ Stanley,
Enumerative Combinatorics, vol.\ 2,
Cambridge Univ.\ Press, Cambridge/New York, 1999.


\bibitem{Y}
A.\ Young,
{\em Quantitative substitutional analysis II},
Proc.\ London Math.\ Soc., Ser.~1, {\bf 35} (1902), 361--397.

\end{thebibliography}
\end{document}